\documentclass[a4paper]{article}
\usepackage{amsmath,amssymb,amsfonts}

\newtheorem{theorem}{Theorem}[section]

\title{The cardinality of the set of real numbers}

\author{Jailton C. Ferreira}

\date{ }
\begin{document}
\maketitle
\pagenumbering{arabic}

\begin{abstract}
\begin{center}
A proof that the set of real numbers is denumerable is given.
\end{center}
\end{abstract}

\section{Introduction} \label{sec-1}

\hspace{22pt}Cantor's proof that the real set is uncountable is now included in the standard mathematics curriculum. The objections and reservations presented by Poincar\'{e}, Kronecker, Weil and Brouwer, for example, are silenced. Apparently we are in the Paradise that Cantor has created. This paper seeks to show that this is not the case.

\hspace{22pt}Each real number of the interval [0, 1] can be
represented by an infinite path in a given binary tree. In Section
\textbf{~\ref{sec-2}} the binary tree is projected on a grid $N
\times N$ and it is shown that the set of the infinite paths
corresponds one-to-one to the set $N$. The Theorems
~\ref{teorema-1} and ~\ref{teorema-2} give the first proof and the
Theorem ~\ref{teorema-3} provides a second proof.

\hspace{22pt}Section \textbf{~\ref{sec-3}} contains the 1891 Cantor's proof and section \textbf{~\ref{sec-4}} examines the the inconsistency of this proof.


\section{The proof of   $\lvert F \lvert  =  \lvert N \lvert$} \label{sec-2}

\begin{theorem} \label{teorema-1}
Let B be a binary tree such that every node has two children and
its depth is equal to $\lvert N \lvert$. Let A be the set of the
nodes of the binary tree B. The cardinality of A is less than or
equal to $\lvert N \lvert$.
\end{theorem}
\textit{Proof:}

\hspace{22pt}Let us now consider the Figure 1. The horizontal
sequence of finite natural numbers, presented in increasing order
from left to right, contains all numbers of $N$; to each natural
number of the sequence corresponds a vertical line. The vertical
sequence of numbers, presented in increasing order from top to
bottom, contains all numbers of $N$; to each natural number of the
sequence corresponds a horizontal line. To each node of the grid
formed by the horizontal and vertical lines corresponds a pair
($m$, $n$), where $m$ belongs to the horizontal sequence of
numbers and $n$ belongs to the vertical sequence of numbers.

\begin{center}
\begin{picture}(180, 180)
\linethickness{0.1mm}

\put (15,15) {\line(1,0){165}} \put (15,30) {\line(1,0){165}} \put
(15,45) {\line(1,0){165}} \put (15,60) {\line(1,0){165}} \put
(15,75) {\line(1,0){165}} \put (15,90) {\line(1,0){165}} \put
(15,105) {\line(1,0){165}} \put (15,120) {\line(1,0){165}} \put
(15,135) {\line(1,0){165}} \put (15,150) {\line(1,0){165}}

\put (30,0) {\line(0,1){165}} \put (45,0) {\line(0,1){165}} \put
(60,0) {\line(0,1){165}} \put (75,0) {\line(0,1){165}} \put (90,0)
{\line(0,1){165}} \put (105,0) {\line(0,1){165}} \put (120,0)
{\line(0,1){165}} \put (135,0) {\line(0,1){165}} \put (150,0)
{\line(0,1){165}} \put (165,0) {\line(0,1){165}}

\put (6,13) {\small 9} \put (6,28) {\small 8} \put (6,43) {\small
7} \put (6,58) {\small 6} \put (6,73) {\small 5} \put (6,88)
{\small 4} \put (6,103) {\small 3} \put (6,118) {\small 2} \put
(6,133) {\small 1} \put (6,148) {\small 0}

\put (29,171) {\small 0} \put (44,171) {\small 1} \put (59,171)
{\small 2} \put (74,171) {\small 3} \put (89,171) {\small 4} \put
(104,171) {\small 5} \put (119,171) {\small 6} \put (134,171)
{\small 7} \put (149,171) {\small 8} \put (165,171) {\small 9}

\end{picture}
\end{center}
\begin{center}
\textrm{Figure 1}
\end{center}

Let us project the binary tree $B$ on the grid of the Figure 1 in
the following way: to the root of the tree corresponds the pair
(0, 0); to the 2 nodes of the level 1 correspond the pairs (0, 1)
and (1, 0); to the 4 nodes of the level 2 correspond the pairs (0,
3), (1, 2), (2, 1) and (3, 0); to the 8 nodes of the level 3
correspond the pairs (0, 7), (1, 6), (2,5), (3, 4), (4, 3), (5,
2), (6, 1) and (7, 0) and to the $2^k$ nodes of the level $k$
correspond the pairs

\begin{equation}\label{dois-1}
(0, 2^{k}-1), (1, 2^{k}-2), (2, 2^{k}-3),  \ldots ,(2^{k}-3, 2),
(2^{k}-2, 1) \hspace{8pt} \textrm{and} \hspace{8pt} (2^{k}-1, 0)
\end{equation}

The Figure 2 shows the binary tree $B$ up to the depth 4. The
Figure 3 shows the projection of the binary tree up to the depth 3
on the grid of the Figure 1.

\begin{center}
\begin{picture}(160, 190)

\put (120,175) {\line(6,1){30}} \put (120,175) {\line(6,-1){30}}
\put (120,155) {\line(6,1){30}} \put (120,155) {\line(6,-1){30}}
\put (120,135) {\line(6,1){30}} \put (120,135) {\line(6,-1){30}}
\put (120,115) {\line(6,1){30}} \put (120,115) {\line(6,-1){30}}
\put (120,95) {\line(6,1){30}} \put (120,95) {\line(6,-1){30}}
\put (120,75) {\line(6,1){30}} \put (120,75) {\line(6,-1){30}}
\put (120,55) {\line(6,1){30}} \put (120,55) {\line(6,-1){30}}
\put (120,35) {\line(6,1){30}} \put (120,35) {\line(6,-1){30}}

\put (80,165) {\line(4,1){40}} \put (80,165) {\line(4,-1){40}}
\put (80,125) {\line(4,1){40}} \put (80,125) {\line(4,-1){40}}
\put (80,85) {\line(4,1){40}} \put (80,85) {\line(4,-1){40}} \put
(80,45) {\line(4,1){40}} \put (80,45) {\line(4,-1){40}}

\put (40,145) {\line(2,1){40}} \put (40,145) {\line(2,-1){40}}
\put (40,65) {\line(2,1){40}} \put (40,65) {\line(2,-1){40}}

\put (0,105) {\line(1,1){40}} \put (0,105) {\line(1,-1){40}}

\put (0,95) {.} \put (40,135) {\scriptsize 0} \put (40,55)
{\scriptsize 1} \put (80,155) {\scriptsize 0} \put (80,115)
{\scriptsize 1} \put (80,75) {\scriptsize 0} \put (80,35)
{\scriptsize 1} \put (120,165) {\scriptsize 0} \put (120,145)
{\scriptsize 1} \put (120,125) {\scriptsize 0} \put (120,105)
{\scriptsize 1} \put (120,85) {\scriptsize 0} \put (120,65)
{\scriptsize 1} \put (120,45) {\scriptsize 0} \put (120,25)
{\scriptsize 1}

\put (40,0) {\small 1} \put (80,0) {\small 2} \put (120,0) {\small
3} \put (150,0) {\small 4}

\end{picture}
\end{center}
\begin{center}
\textrm{Figure 2}
\end{center}

\begin{center}
\begin{picture}(180, 180)
\linethickness{0.1mm}

\put (15,15) {\line(1,0){165}} \put (15,30) {\line(1,0){165}} \put
(15,45) {\line(1,0){165}} \put (15,60) {\line(1,0){165}} \put
(15,75) {\line(1,0){165}} \put (15,90) {\line(1,0){165}} \put
(15,105) {\line(1,0){165}} \put (15,120) {\line(1,0){165}} \put
(15,135) {\line(1,0){165}} \put (15,150) {\line(1,0){165}}

\put (30,0) {\line(0,1){165}} \put (45,0) {\line(0,1){165}} \put
(60,0) {\line(0,1){165}} \put (75,0) {\line(0,1){165}} \put (90,0)
{\line(0,1){165}} \put (105,0) {\line(0,1){165}} \put (120,0)
{\line(0,1){165}} \put (135,0) {\line(0,1){165}} \put (150,0)
{\line(0,1){165}} \put (165,0) {\line(0,1){165}}

\put (6,13) {\small 9} \put (6,28) {\small 8} \put (6,43) {\small
7} \put (6,58) {\small 6} \put (6,73) {\small 5} \put (6,88)
{\small 4} \put (6,103) {\small 3} \put (6,118) {\small 2} \put
(6,133) {\small 1} \put (6,148) {\small 0}

\put (29,171) {\small 0} \put (44,171) {\small 1} \put (59,171)
{\small 2} \put (74,171) {\small 3} \put (89,171) {\small 4} \put
(104,171) {\small 5} \put (119,171) {\small 6} \put (134,171)
{\small 7} \put (149,171) {\small 8} \put (165,171) {\small 9}

\thicklines \put (30,150){\line(1,0){105}} \put
(30,150){\line(0,-1){105}} \put (45,150){\line(1,-1){45}} \put
(30,135){\line(1,-1){45}} \put (75,150){\line(3,-1){45}} \put
(60,135){\line(3,-1){45}} \put (45,120){\line(1,-3){15}} \put
(30,105){\line(1,-3){15}}

\end{picture}
\end{center}
\begin{center}
\textrm{Figure 3}
\end{center}

\hspace{22pt}It is known that the bijection  $f:N \rightarrow  N
\times N$, where $N$ is the set of finite natural numbers, can be
defined applying the diagonal method.

\begin{center}
\begin{picture}(180, 180)
\linethickness{0.1mm}

\put (15,15) {\line(1,0){165}} \put (15,30) {\line(1,0){165}} \put
(15,45) {\line(1,0){165}} \put (15,60) {\line(1,0){165}} \put
(15,75) {\line(1,0){165}} \put (15,90) {\line(1,0){165}} \put
(15,105) {\line(1,0){165}} \put (15,120) {\line(1,0){165}} \put
(15,135) {\line(1,0){165}} \put (15,150) {\line(1,0){165}}

\put (30,0) {\line(0,1){165}} \put (45,0) {\line(0,1){165}} \put
(60,0) {\line(0,1){165}} \put (75,0) {\line(0,1){165}} \put (90,0)
{\line(0,1){165}} \put (105,0) {\line(0,1){165}} \put (120,0)
{\line(0,1){165}} \put (135,0) {\line(0,1){165}} \put (150,0)
{\line(0,1){165}} \put (165,0) {\line(0,1){165}}

\put (6,13) {\small 9} \put (6,28) {\small 8} \put (6,43) {\small
7} \put (6,58) {\small 6} \put (6,73) {\small 5} \put (6,88)
{\small 4} \put (6,103) {\small 3} \put (6,118) {\small 2} \put
(6,133) {\small 1} \put (6,148) {\small 0}

\put (29,171) {\small 0} \put (44,171) {\small 1} \put (59,171)
{\small 2} \put (74,171) {\small 3} \put (89,171) {\small 4} \put
(104,171) {\small 5} \put (119,171) {\small 6} \put (134,171)
{\small 7} \put (149,171) {\small 8} \put (165,171) {\small 9}

\thicklines \put (30,150) {\line(1,0){15}} \put (45,150)
{\line(-1,-1){15}} \put (30,135) {\line(0,-1){15}} \put (30,120)
{\line(1,1){30}} \put (60,150) {\line(1,0){15}} \put (75,150)
{\line(-1,-1){45}} \put (30,105) {\line(0,-1){15}} \put (30,90)
{\line(1,1){60}} \put (90,150) {\line(1,0){15}} \put (105,150)
{\line(-1,-1){75}} \put (30,75) {\line(0,-1){15}} \put (30,60)
{\line(1,1){90}} \put (30,90) {\line(1,1){60}} \put (120,150)
{\line(1,0){15}} \put (135,150) {\line(-1,-1){105}} \put (30,45)
{\line(0,-1){15}} \put (30,30) {\vector(1,1){90}}

\end{picture}
\end{center}
\begin{center}
\textrm{Figure 4}
\end{center}

The Figure 4 shows how starting from the pair (0, 0) and following
the thick line we can establish the one-to-one correspondence
between the set of the pairs ($m$, $n$) and $N$, this is,

\begin{equation}\label{dois-2}
\begin{tabular}{c c c c c c c c}
(0,0) & (1,0) & (0,1) & (0,2) & (1,1) & (2,0) & (3,0) & \ldots \\
 $\updownarrow$ & $\updownarrow$ & $\updownarrow$ & $\updownarrow$ & $\updownarrow$ & $\updownarrow$ & $\updownarrow$ &  \\
0 & 1 & 2 & 3 & 4 & 5 & 6 & \ldots \\
\end{tabular}
\end{equation}

Considering that to any node of the binary tree corresponds a node
of the grid and that there are nodes in the grid without
corresponding nodes in the binary tree, we conclude that

\begin{equation}\label{dois-3}
\lvert A \lvert  \le \lvert N \lvert
\end{equation}

\begin{theorem} \label{teorema-2}
The cardinality of the set of real numbers of the interval [0, 1]
is equal to $\lvert N \lvert$.
\end{theorem}
\textit{Proof:}

\hspace{22pt}Let us denote by $F$ the set of real numbers of the
interval [0, 1]. Any element of $F$ can be represented in the
binary system by

\begin{equation}\label{dois-4}
f_1 \times 2^{-1} + f_2 \times 2^{-2} + f_3 \times 2^{-3} + f_4
\times 2^{-4} + \ldots
\end{equation}

where
\begin{equation}\label{dois-5}
f_i \in \{ 0, 1 \} \qquad \textrm{and} \qquad i \in N \qquad
\textrm{and} \qquad i \ne 0
\end{equation}

The representation \eqref{dois-4} can be simplified to

\begin{equation}\label{dois-6}
.f_1 f_2 f_3 f_4 f_5 \ldots
\end{equation}

where the first character of the sequence is the point. Each $f_i$
of \eqref{dois-6} is substituted by 0 or 1 to represent a given
number. The $i$-th 0 or 1 on the right of ``." corresponds to
$f_i$. The set $F$ can also be represented by a binary tree where
each node has two children. Each infinite path on the binary tree,
$.q_{\tau 1} q_{\tau 2} \ldots  q_{\tau i} \ldots $ with depth
equal to $\lvert N \lvert$ , represents an element of $F$.

\hspace{22pt} Let us denote by $B$ the binary tree above
mentioned; by $B_{i}$ the binary tree of depth $i$; by $P$ the set
of all paths starting from the root of $B$; by $P_i$ the set of
all paths from the root of $B_i$ to its leaves. Each element of
$P_i$ or $P$ is the set of the nodes (except the root) that form
the path. Let the set

\begin{equation}\label{dois-7}
\{P_{1}, P_{2}, \ldots , P_{i}, \ldots \} \qquad \textrm{such
that} \qquad i \in N
\end{equation}

The set $\{P_{1}, P_{2}, \ldots , P_{i} \}$ contains all paths
from the root of the tree $B_{i}$ and the endings of the paths that
exist in $\{P_{1}, P_{2}, \ldots , P_{i} \}$ are all the nodes of
$B_i$, except the root. This indicates that the set in
\eqref{dois-7} - which corresponds to all the tree $B$ - contains
all paths of $B$.

\hspace{22pt} Let us consider the set of all infinite paths from
the root of $B$

\begin{equation}\label{dois-8}
Q=\{q_{\alpha}, q_{\beta}, q_{\gamma}, \ldots, q_{\tau}, \ldots\}
\end{equation}

The infinite path $q_{\tau}$ can be represented by the set

\begin{equation}\label{dois-9}
\{q_{\tau 1}, q_{\tau 2}, \ldots , q_{\tau i}, \ldots \}
\end{equation}

where $q_{\tau i}$ is the $i$-th node of $q_{\tau}$ after the
root. Let

\begin{equation}\label{dois-10}
Q_{\tau} = \{ \{q_{\tau 1} \}, \{q_{\tau 1} , q_{\tau 2}\}, \ldots
, \{q_{\tau 1}, q_{\tau 2}, \ldots , q_{\tau i} \}, \ldots \}
\end{equation}

The union of all elements of $Q_{\tau}$ is equal to $q_{\tau}$,
that is,

\begin{equation}\label{dois-11}
q_{\tau} = \{q_{\tau 1} \} \hspace{4pt} \cup \hspace{4pt}
\{q_{\tau 1} , q_{\tau 2}\} \hspace{4pt} \cup \hspace{4pt} \ldots
\hspace{4pt} \cup \hspace{4pt} \{q_{\tau 1}, q_{\tau 2}, \ldots ,
q_{\tau i} \} \hspace{4pt} \cup \hspace{4pt} \ldots
\end{equation}

Considering \eqref{dois-11} we can obtain from \eqref{dois-8}

\begin{eqnarray}\label{dois-12}
Q = \{ \{q_{\alpha 1} \} \hspace{4pt} \cup \hspace{4pt}
\{q_{\alpha 1} , q_{\alpha 2}\} \hspace{4pt} \cup \hspace{4pt}
\ldots \hspace{4pt} \cup \hspace{4pt} \{q_{\alpha 1}, q_{\alpha
2}, \ldots , q_{\alpha i} \}
\hspace{4pt} \cup \hspace{4pt} \ldots , \nonumber\\
\{q_{\beta 1} \} \hspace{4pt} \cup \hspace{4pt} \{q_{\beta 1} ,
q_{\beta 2}\} \hspace{4pt} \cup \hspace{4pt} \ldots \hspace{4pt}
\cup \hspace{4pt} \{q_{\beta 1}, q_{\beta 2}, \ldots , q_{\beta i}
\} \hspace{4pt} \cup \hspace{4pt} \ldots , \nonumber\\ \textit{.}
\ldots \ldots \ldots \ldots \ldots \ldots \ldots \ldots \ldots
\ldots \ldots \ldots \ldots \ldots \ldots \ldots \ldots \ldots \textit{.} \nonumber\\
\{q_{\tau 1} \} \hspace{4pt} \cup \hspace{4pt} \{q_{\tau 1} ,
q_{\tau 2}\} \hspace{4pt} \cup \hspace{4pt} \ldots \hspace{4pt}
\cup \hspace{4pt} \{q_{\tau 1}, q_{\tau 2}, \ldots ,
q_{\tau i} \} \hspace{4pt} \cup \hspace{4pt} \ldots , \nonumber\\
\ldots \ldots \ldots \ldots \ldots \ldots \ldots \ldots \ldots
\ldots \ldots \ldots \ldots \ldots \ldots \ldots \ldots
\ldots \textit{.} \} \nonumber\\
\end{eqnarray}

From \eqref{dois-12} we have

\begin{eqnarray}\label{dois-13}
\lvert Q \lvert \leq \lvert \{q_{\alpha 1} \} \hspace{4pt} \cup
\hspace{4pt} \{q_{\alpha 1} , q_{\alpha 2}\} \hspace{4pt} \cup
\hspace{4pt} \ldots \hspace{4pt} \cup \hspace{4pt} \{q_{\alpha 1},
q_{\alpha 2}, \ldots , q_{\alpha i} \}
\hspace{4pt} \cup \hspace{4pt} \ldots \nonumber\\
\cup \hspace{4pt} \{q_{\beta 1} \} \hspace{4pt} \cup \hspace{4pt}
\{q_{\beta 1} , q_{\beta 2}\} \hspace{4pt} \cup \hspace{4pt}
\ldots \hspace{4pt} \cup \hspace{4pt} \{q_{\beta 1}, q_{\beta 2},
\ldots , q_{\beta i} \} \hspace{4pt} \cup \hspace{4pt} \ldots \nonumber\\
\ldots \ldots \ldots \ldots \ldots \ldots \ldots \ldots \ldots
\ldots \ldots \ldots \ldots \ldots \ldots \ldots \ldots \ldots
\ldots \nonumber\\
\cup \hspace{4pt} \{q_{\tau 1} \} \hspace{4pt} \cup \hspace{4pt}
\{q_{\tau 1} , q_{\tau 2}\} \hspace{4pt} \cup \hspace{4pt} \ldots
\hspace{4pt} \cup \hspace{4pt} \{q_{\tau 1}, q_{\tau 2}, \ldots ,
q_{\tau i} \} \hspace{4pt} \cup \hspace{4pt} \ldots \nonumber\\
\ldots \ldots \ldots \ldots \ldots \ldots \ldots \ldots \ldots
\ldots \ldots \ldots \ldots \ldots \ldots \ldots \ldots \ldots
\ldots \lvert \nonumber\\
\end{eqnarray}

Since the sets on the right side of \eqref{dois-13} with
cardinality equal to $i$ are the elements of $P_i$, that is,

\begin{equation}\label{dois-14}
P_{i} = \{q_{\alpha 1}, q_{\alpha 2}, \ldots , q_{\alpha i} \}
\hspace{4pt} \cup \hspace{4pt} \{q_{\beta 1}, q_{\beta 2}, \ldots
, q_{\beta i} \} \hspace{4pt} \cup \hspace{4pt} \{q_{\tau 1},
q_{\tau 2}, \ldots , q_{\tau i} \} \hspace{4pt} \cup \ldots
\end{equation}

we conclude that

\begin{equation}\label{dois-15}
\lvert Q \lvert \leq \lvert P_{1} \lvert \hspace{4pt} +
\hspace{4pt} \lvert P_{2} \lvert \hspace{4pt} + \hspace{4pt}
\lvert P_{3} \lvert \hspace{4pt} + \hspace{4pt} \ldots
\hspace{4pt} + \lvert P_{i} \lvert \hspace{4pt} + \ldots
\end{equation}

\hspace{22pt} The set of the nodes of $B$ has cardinality less
than or equal to $\lvert N \lvert$ (Theorem 2.1). Since the number
of nodes in the level $i$ of the tree $B_i$ is equal to $\lvert
P_{i} \lvert$, we have

\begin{equation}\label{dois-16}
\lvert N \lvert \ge \lvert P_{1} \lvert \hspace{4pt} +
\hspace{4pt} \lvert P_{2} \lvert \hspace{4pt} + \hspace{4pt}
\lvert P_{3} \lvert \hspace{4pt} + \hspace{4pt} \ldots
\hspace{4pt} + \lvert P_{i} \lvert \hspace{4pt} + \ldots
\end{equation}

Therefore

\begin{equation}\label{dois-17}
\lvert Q \lvert \leq \lvert N \lvert
\end{equation} \

Considering the bijection $f:F \rightarrow Q$ and $\lvert F \lvert
\geq \lvert N \lvert$, we conclude

\begin{equation}\label{dois-18}
\lvert F \lvert = \lvert N \lvert
\end{equation}

\begin{theorem} \label{teorema-3}
The cardinal number of the set of infinite paths of the binary
tree B is equal to $ \lvert N \lvert $.
\end{theorem}
\textit{Proof:}

\hspace{22pt}Be the infinite path of $Q$

\begin{equation}\label{dois-19}
q_{\tau 0},q_{\tau 1}, q_{\tau 2}, \ldots , q_{\tau i}, \ldots
\end{equation}

where $q_{\tau i} = 1$ for all $i \in N$ with $i > 0$, projected in the grid as
shows the Figure 3. In the coordinates of the grid, the path
\eqref{dois-19} corresponds to the pairs

\begin{equation}\label{dois-20}
(0, 0), (0, 1), (0, 3), (0, 7), (0, 15), \ldots , (0, 2^{k}-1),
\ldots
\end{equation}

Let us notice that in (0, 0) all infinite paths begin. By the pair
(0, 1) a part of all infinite paths passes, by (0, 3) a part of
the infinite paths that passed by (0, 1) passes, by (0, 7) a part
of the infinite paths that passed by (0, 3) passes and so forth.
When an entire path is accomplished, the path exists. In the
examined case, \eqref{dois-19} is the accomplished path.

\hspace{22pt}Each pair $(0, 2^{k}-1)$ of the sequence
\eqref{dois-20} belongs to the set of pairs given by
\eqref{dois-21}

\begin{equation}\label{dois-21}
(0, 2^{k}-1), (1, 2^{k}-2), (2, 2^{k}-3),  \ldots , (2^{k}-3, 2),
(2^{k}-2, 1) \hspace{8pt} \textrm{and} \hspace{8pt} (2^{k}-1, 0)
\end{equation}

Let us denote by $G_k$ the set of pairs given by \eqref{dois-21}
for $k$. For any $k$ the cardinal number of the set of infinite
paths that passes by $(0, 2^{k}-1)$ is equal to the cardinal
number of the set of infinite paths that pass by any other pair of
\eqref{dois-21}. For any $k$ the distance from (0, 0) to $(0,
2^{k}-1)$ in the grid is equal to $\lvert G_k \lvert$. When we
examine the pairs of \eqref{dois-20}, from left to right, and we
accomplished the path \eqref{dois-19}, $G_k$ becomes the set of
the infinite path whose cardinality is equal to the cardinality of
the set of the nodes of the path \eqref{dois-19}.

\hspace{22pt}The cardinality of the set of the infinite paths is
$\lvert F \lvert$ and the cardinality of the set of nodes of any
infinite path is $\lvert N \lvert$. Therefore

\begin{equation}\label{dois-22}
\lvert F \lvert = \lvert N \lvert
\end{equation}


\section{The 1891 proof} \label{sec-3}

\begin{quote}

\textit{Theorem.} \hspace{8pt} Let $F$ be the set of real numbers of the
interval [0, 1].  The set $F$ is not countable.

\textit{Proof.}\hspace{8pt} We assume that the set $F$ is
countable. This means, by the definition of countable sets, that
$F$ is finite or denumerable. Let us notice that $ \lvert F \lvert
$ is not less than $ \lvert N \lvert $ where $N$ the set of natural numbers. Be

\begin{center}
$F = \{a_1, a_2, a_3, \ldots \}$
\end{center}

where the cardinality of $ \{a_0, a_1, a_2, \ldots \} $ is $
\lvert N \lvert $. We can write their decimal expansions as
follows:

\begin{equation}\label{3-1}
\begin{tabular}{c c c c c c c}
$a_1$ & = & 0. & $d_{1,1}$ & $d_{1,2}$ & $d_{1,3}$ & \ldots \\
$a_2$ & = & 0. & $d_{2,1}$ & $d_{2,2}$ & $d_{2,3}$ & \ldots \\
$a_3$ & = & 0. & $d_{3,1}$ & $d_{3,2}$ & $d_{3,3}$ & \ldots \\
\ldots & & \ldots
\end{tabular}
\end{equation}

where the $d$'s are binary characters 0 and 1. Now we define the
number

\begin{center}
$x = 0.d_1 d_2 d_3 \ldots $
\end{center}

by selecting $d_1 \ne d_{1,1}$,  $d_2 \ne d_{2,2}$,  $d_3 \ne
d_{3,3}$, \ldots . This gives a number not in the set $ \{a_1,
a_2, a_3, \ldots \} $, but $x \in F$. Therefore, $F$ is not
countable ~\cite{Cantor}.
\end{quote}


\section{The inconsistency of the proof} \label{sec-4}

\hspace{22pt} When \eqref{3-1} is represented as \eqref{4-1} bellow, the
proof can be called the \textit{written list form} of the Cantor's
argument of 1891.

\begin{equation}\label{4-1}
\begin{tabular}{c c c c c c c}
$1$ & $\longleftrightarrow$ & 0. & $d_{1,1}$ & $d_{1,2}$ & $d_{1,3}$ & \ldots \\
$2$ & $\longleftrightarrow$ & 0. & $d_{2,1}$ & $d_{2,2}$ & $d_{2,3}$ & \ldots \\
$3$ & $\longleftrightarrow$ & 0. & $d_{3,1}$ & $d_{3,2}$ & $d_{3,3}$ & \ldots \\
\ldots & & \ldots
\end{tabular}
\end{equation}

Let be the list

\begin{equation}\label{4-2}
\begin{tabular}{c c c c}
$n \in N$ & $\longleftrightarrow$ & $x \in F$ \\
\\
$0$ & $\longleftrightarrow$ & 0.000000 & \ldots \\
$1$ & $\longleftrightarrow$ & 0.100000 & \ldots \\
$2$ & $\longleftrightarrow$ & 0.010000 & \ldots \\
$3$ & $\longleftrightarrow$ & 0.110000 & \ldots \\
$4$ & $\longleftrightarrow$ & 0.001000 & \ldots \\
$5$ & $\longleftrightarrow$ & 0.101000 & \ldots \\
$6$ & $\longleftrightarrow$ & 0.011000 & \ldots \\
$7$ & $\longleftrightarrow$ & 0.111000 & \ldots \\
$8$ & $\longleftrightarrow$ & 0.000100 & \ldots \\
$9$ & $\longleftrightarrow$ & 0.100100 & \ldots \\
$10$ & $\longleftrightarrow$ & 0.010100 & \ldots \\
$11$ & $\longleftrightarrow$ & 0.110100 & \ldots \\
$12$ & $\longleftrightarrow$ & 0.001100 & \ldots \\
$13$ & $\longleftrightarrow$ & 0.101100 & \ldots \\
$14$ & $\longleftrightarrow$ & 0.011100 & \ldots \\
$15$ & $\longleftrightarrow$ & 0.111100 & \ldots \\
$16$ & $\longleftrightarrow$ & 0.000010 & \ldots \\
\ldots & & \ldots
\end{tabular}
\end{equation}

Let us notice that (i) $d_{1,1}$, $d_{3,1}$, $d_{5,1}$, $d_{7,1}$, $d_{9,1}$, $d_{11,1}$, $d_{13,1}$, $d_{15,1}$, $d_{17,1}$ ... are equal to 0 and  $d_{2,1}$, $d_{4,1}$, $d_{6,1}$, $d_{8,1}$, $d_{10,1}$, $d_{12,1}$, $d_{14,1}$, $d_{16,1}$, $d_{18,1}$ ... are equal to 1, (ii) $d_{1,2}$, $d_{2,2}$, $d_{5,2}$, $d_{6,2}$, $d_{9,2}$, $d_{10,2}$, $d_{13,2}$, $d_{14,2}$, $d_{17,2}$ ... are equal to 0 and  $d_{3,2}$, $d_{4,2}$, $d_{7,2}$, $d_{8,2}$, $d_{11,2}$, $d_{12,2}$, $d_{15,2}$, $d_{16,2}$, $d_{19,2}$ ... are equal to 1,  (iii) $d_{1,3}$, $d_{2,3}$, $d_{3,3}$, $d_{4,3}$, $d_{9,3}$, $d_{10,3}$, $d_{11,3}$, $d_{12,3}$, $d_{17,3}$ ... are equal to 0 and  $d_{5,3}$, $d_{6,3}$, $d_{7,3}$, $d_{8,3}$, $d_{13,3}$, $d_{14,3}$, $d_{15,3}$, $d_{16,3}$, ... are equal to 1, (iv) ...  and so on ...

\begin{equation}\label{4-3}
\begin{tabular}{c c cc c c}
 $d_{1,1}$ & $d_{1,2}$ & $d_{1,3}$ & \ldots \\
 $d_{2,1}$ & $d_{2,2}$ & $d_{2,3}$ & \ldots \\
$d_{3,1}$ & $d_{3,2}$ & $d_{3,3}$ & \ldots \\
\ldots
\end{tabular}
\end{equation}

More generally, considering the matrix \eqref{4-3} the procedure to obtain \eqref{4-2} is:  the first column of the matrix is filled from top to
bottom with a succession of the pattern 01. The second column is
filled with a succession of the pattern 0011. The $n$-th column is filled from top to bottom with a succession of
the pattern constituted of $2^{n}$ 0's followed by $2^{n}$ 1's.

\hspace{22pt} When the number of rows in \eqref{4-2} goes to infinity we have that a $x$ equal to

\begin{equation}\label{4-4}
lim_{n \to \infty} (2^{-1}+2^{-2}+2^{-3}+ ... +2^{-n})= 0.11111111...
\end{equation}

belongs to \eqref{4-2}. There is no way out of this. All the reals in the interval [0,1] are included in the list. However applying the Cantor's argument to \eqref{4-2} we found that 0.11111111... does not belongs to \eqref{4-2}.

\end{document}